\documentclass{paper}

\usepackage{amssymb}
\usepackage{amscd}
\usepackage{amsfonts}
\usepackage{amssymb}
\usepackage{cite}
\usepackage{amsmath}
\usepackage{epsfig}
\usepackage{graphicx, subfigure}
\usepackage{color, graphics}
\usepackage{pdfsync}
\usepackage{hyperref}
\usepackage[all,cmtip]{xy}
\usepackage{setspace}

\newtheorem{teo}{\textbf{Theorem}}
\newtheorem{prp}{\textbf{Proposition}}
\newtheorem{lem}{\textbf{Lemma}}

\newtheorem{rem}{\textbf{Remark}}

\newcommand{\e}{\mathbf{e}}
\newcommand{\p}{\mathbf{p}}

\newcommand{\R}{\mathbb{ R}}

\DeclareMathOperator{\Span}{span}

\DeclareMathOperator{\tr}{tr}
\DeclareMathOperator{\diag}{diag}

\numberwithin{equation}{section}

\allowdisplaybreaks

\begin{document}
\title{Contact magnetic geodesic and sub-Riemannian flows on $V_{n,2}$ and integrable cases of a heavy rigid body with a gyrostat}

\author{Bo\v zidar Jovanovi\' c}

\maketitle

\noindent{\small Mathematical Institute SANU, Kneza Mihaila 36, Belgrade, Serbia\footnote{{\sc email:} bozaj@mi.sanu.ac.rs}}

\begin{abstract}
We prove the integrability of magnetic geodesic flows of $SO(n)$--invariant Riemannian metrics on the rank two Stefel variety $V_{n,2}$ with respect to the magnetic field $\eta\, d\alpha$, where $\alpha$ is the standard contact form on $V_{n,2}$ and $\eta$ is a real parameter.
Also, we prove the integrability of magnetic sub-Riemannian geodesic flows for $SO(n)$-invariant sub-Riemannian structures on $V_{n,2}$. All statements in the limit $\eta=0$ imply the integrability of the problems without the influence of the magnetic field. We also consider integrable pendulum-type natural mechanical systems with the kinetic energy defined by $SO(n)\times SO(2)$--invariant Riemannian metrics. For $n=3$, using the isomorphism $V_{3,2}\cong SO(3)$, the obtained integrable magnetic models reduce to
integrable cases of a motion of a heavy rigid body with a gyrostat around a fixed point: Zhukovskiy--Volterra gyrostat, the Lagrange top with a gyrostat, and the Kowalevski top with a gyrostat. As a by-product we obtain the Lax presentations for the Lagrange gyrostat and the Kowalevski gyrostat in the fixed reference frame (dual Lax representations).\footnote{{\sc msc:} 37J35;  53D25; 53C17; 70E40, 70G65. {\sc keywords:} {magnetic geodesic and sub-Riemannian flows;  Liouville and noncommutative integrability; contact structure; Zhukovskiy--Volterra gyrostat, the Lagrange top; the Kowalevski top}.}
\end{abstract}

\section{Introduction}

The rank two Stiefel variety $V_{n,2}=SO(n)/SO(n-2)$ is the variety of ordered sets of two orthogonal unit  vectors $\e_1,\e_2$ in the
Euclidean space $(\R^n,\langle\cdot,\cdot\rangle)$:
\begin{equation}\label{e1e2}
V_{n,2}: \qquad  \langle \e_1,\e_1\rangle=1, \qquad \langle
\e_2,\e_2\rangle=1, \qquad \langle \e_1,\e_2 \rangle=0.
\end{equation}

The tangent bundle $TV_{n,2}$ is a ($4n-6$)--dimensional
subspace of $\R^{4n}(\e_1,\e_2,\xi_1,\xi_2)$ defined by additional constraints
\begin{equation}\label{e-dot}
\langle \e_1,\xi_1\rangle=0, \qquad \langle \e_2,\xi_2\rangle=0, \qquad \langle \e_1,\xi_2 \rangle+\langle \xi_1,\e_2\rangle=0.
\end{equation}

The Stiefel variety $V_{n,2}$ can be seen also as the unit sphere bundle $T_1 S^{n-1}$ with respect to the standard round sphere metric
\[
T_1 S^{n-1}=\{\e_2\in T_{\e_1} S^{n-1}\, \vert \langle \e_2,\e_2\rangle=1, \, \e_1\in S^{n-1}\}, \quad S^{n-1}=\{\e_1\in \R^n\,\vert\, \langle \e_1,\e_1\rangle=1\}.
\]
Therefore, it carries the \emph{standard contact form}, the restriction of the Liouville 1-form from $TS^{n-1}\cong T^*S^{n-1}$  to the unit tangent bundle
\[
\alpha=-\e_2 d\e_1\vert_{V_{n,2}}=-\sum_{i=1}^n e_2^i de_1^i\vert_{V_{n,2}}.
\]

Let $\mathcal H=\ker\alpha\subset TV_{n,2}$ be the \emph{standard contact distribution}. We have
\[
\xi=(\xi_1,\xi_2)\in  \mathcal H\vert_{(\e_1,\e_2)}\subset T_{(\e_1,\e_2)} V_{n,2} \quad \Longleftrightarrow \quad  \sum_{i=1}^n e_2^i de_1^i(\xi)=\langle \e_2,\xi_1 \rangle=0.
\]
Taking into account the constraints \eqref{e-dot},  we obtain that $\mathcal H$ is given by
\[
\mathcal H =\big\{(\e_1,\e_2,\xi_1,\xi_2)\in V_{n,2}\times \R^{2n}\,\vert\,\langle \e_1,\xi_1\rangle=0,\,
\langle \e_2,\xi_2\rangle=0, \, \langle \e_1,\xi_2 \rangle=0, \, \langle \xi_1,\e_2\rangle=0\big\}
\]
that is,
$
\mathcal H\vert_{(\e_1,\e_2)}=(\Span\{\e_1,\e_2\}^\perp,\Span\{\e_1,\e_2\}^\perp).
$

The fact that  $(V_{n,2},\alpha)$ is a contact manifold is equivalent to the non-degeneracy of the closed two form
\[
\omega_{mag}=d\alpha=d\e_1 \wedge d\e_2\vert_{V_{n,2}}=\sum_{i=1}^n de_1^i \wedge de_2^j\vert_{V_{n,2}}
\]
restricted to $\mathcal H$, or, to the condition $\alpha\wedge \omega_{mag}^{n-2}\ne 0$ (see \cite{LM}).  We refer to $\omega_{mag}$ as the \emph{standard contact magnetic form} on $V_{n,2}$.

Thus, we can study the following four natural problems on the rank two Stiefel variety:

\begin{itemize}

\item{}  Magnetic geodesic flows with respect to the magnetic force defined by $\eta\,\omega_{mag}$. Here $\eta$ is a real parameter representing the strength of the magnetic field, and for $\eta=0$ we have the usual geodesic flows.

\item{} Sub-Riemannian magnetic geodesic flows, with the sub-Riemannian structures defined on $\mathcal H$ and other $SO(n)$--invariant bracket generating distributions.

\item{} Natural mechanical systems with influence of the magnetic field defined by $\eta\,\omega_{mag}$.

\item{} Contact Hamiltonian flows on $(V_{n,2},\alpha)$.

\end{itemize}

In this paper we prove the integrability of magnetic geodesic flows of $SO(n)$--invariant Riemannian metrics (Theorem \ref{glavna}, Section \ref{sec3I}).
We also describe $SO(n)$--invariant sub-Riemannian structures on the contact distribution $\mathcal H$ (Section \ref{secSRH}), as well as on other invariant bracket generating distributions (Section \ref{secSRD}).  We consider normal sub-Riemannian geodesics (see \cite{M2002, ABB}) and prove integrability of the corresponding magnetic geodesic flows (Theorem \ref{glavna2}, Section \ref{sec3II}). In particular, by the restriction of the Euclidean metric from the ambient space $\R^{2n}(\e_1,\e_2)$ to the contact distribution $\mathcal H$ we obtain the standard sub-Riemannian structure $ds^2_{\mathcal H}$ with integrable magnetic flow.
All statements in the limit $\eta=0$ imply the integrability of the problems without the influence of the magnetic field.

In \cite{FeJo_MZ_2012, FeJo_JPA_2012} we  considered potential systems with $SO(n)\times SO(2)$--invariant Riemannian metrics and proved the integrability for the case of the Neumann potential and a pendulum-type potential, the pull-back of the pendulum potential on the oriented Grassmannian variety $G^+_{n,2}$ realized as an adjoint orbit in $so(n)$. The later result is recalled in Section \ref{sec4I} (Theorem \ref{Gn2-klatno}). There is another natural choice for potential function, induced from two independent pendulum systems on $S^{n-1}(\e_1)$ and $S^{n-2}(\e_2)$ (Section \ref{sec4II}). We obtain the Lax representation of the problem for a special $SO(n)\times SO(2)$--invariant metric (Theorem \ref{Vn2-klatno}) which implies the complete integrability of the system (Theorem \ref{main}).

For $n=3$, using the isomorphism $V_{3,2}\cong SO(3)$, the obtained integrable models reduce to integrable cases of a motion of a heavy rigid body with a gyrostat around a fixed point: Zhukovskiy--Volterra gyrostat (Proposition \ref{stav1}, Section \ref{sec3III}), the Lagrange top with a gyrostat (Proposition \ref{stav2}, Section \ref{sec4III}), and the Kowalevski top with a gyrostat (Proposition \ref{stav3}, Section \ref{sec4III}). As a by-product we obtain the Lax presentations for the Lagrange and the Kowalevska tops with gyrostat in the fixed reference frames (dual Lax representations). For the sake of completeness of the exposition, in Section \ref{sec2} we recall the problem of integrability of a motion of a heavy rigid body with a gyrostat.

The study of contact Hamiltonian flows on $(V_{n,2},\alpha)$ is beyond the scope of this paper (for a notion of complete integrability of contact systems, e.g.,  see \cite{JoMatSb} and references therein).

\section{Heavy rigid body with a gyrostat}\label{sec2}
The integrability of a rigid body motion is one of the central problems in classical mechanics.
We recall on the formulation of the problem with addition of a gyroscope. Let us consider a rigid body motion around a fixed point with a gyrostat of constant angular momentum $\mathbf L=(L_1,L_2,L_3)$ in the body frame. Let $I$ be the operator of inertia of the body with a gyrostat (slightly modified rigid body operator, e.g., see \cite{Zh2,  DGJ2023}) and let $\Omega$ be the angular velocity of the body in the moving frame. Then the total angular momentum of the system body+gyrostat is $\mathbf K=\mathbf M+\mathbf L$, where
$\mathbf M=I\Omega$.

Let $\gamma$ be the unit vector fixed in the space in the direction of a homogeneous gravitational field and $\Gamma$ the same vector considered in the frame attached to the body. The Euler-Poisson equation for the motion in the body frame are given by
\begin{equation}\label{MEP}
\begin{aligned}
& \dot{\mathbf M}=\big(\mathbf M+\mathbf L\big) \times \Omega+\Gamma\times \chi, \\
& \dot\Gamma=\Gamma\times\Omega, \quad \Omega=I^{-1}\mathbf M,
\end{aligned}
\end{equation}
where $\chi=(\chi_1,\chi_2,\chi_3)$ is the position of the center of mass of the system body+gyrostat, multiplied by the mass of the system $m$  and the gravitational constant $g$.

Note that the equations \eqref{MEP} are Hamiltonian
\[
\dot M_i=\{M_i,H\}_\mathbf L, \qquad \dot\Gamma_i=\{\Gamma_i,H\}_\mathbf L, \qquad H=\frac12\langle I^{-1}M,M\rangle+\langle \chi,\Gamma\rangle.
 \]
with respect to the magnetic Poisson bracket on $\R^6(\mathbf M,\Gamma)$ obtained from the standard one by the argument translation
\begin{equation}\label{MPB}
\{F,G\}_\mathbf L\vert_{(\mathbf M,\Gamma)}=
-\langle \mathbf M+\mathbf L,\frac{\partial F}{\partial \mathbf M}\times \frac{\partial G}{\partial \mathbf M}\rangle-
\langle \Gamma,\frac{\partial F}{\partial \mathbf M}\times \frac{\partial G}{\partial \Gamma}+\frac{\partial F}{\partial \Gamma}\times \frac{\partial G}{\partial \mathbf M}\rangle.
\end{equation}

It is also convenient to write them with the Hamiltonian having a linear term.
Set $\mathbf K=\mathbf M+\mathbf L$. Then the system \eqref{MEP} is equivalent to the system
\begin{equation}\label{EP}
\begin{aligned}
& \dot{\mathbf K}=\mathbf K\times \Omega+\Gamma\times \chi, \\
& \dot\Gamma=\Gamma\times\Omega, \quad \Omega=I^{-1}\big(\mathbf K-\mathbf L\big),
\end{aligned}
\end{equation}
that is Hamiltonian with respect to the standard Poisson brackets $\{\cdot,\cdot\}_0$ on $\R^6(\mathbf K,\Gamma)$ (set $\mathbf M\mapsto \mathbf K$ and $\mathbf L\mapsto 0$ in \eqref{MPB}):
\[
\dot K_i=\{K_i,H_1\}_0, \qquad \dot\Gamma_i=\{\Gamma_i,H_1\}_0, \qquad H_1=\frac12\langle\mathbf K,I^{-1}\mathbf K\rangle-\langle \mathbf K,I^{-1}\mathbf L\rangle+\langle\Gamma,\chi\rangle.
\]

As in the case without the gyrostat,  in addition to the Hamiltonian function, geometric integral $\langle \Gamma,\Gamma\rangle=1$, and the area integral $\langle \mathbf M+\mathbf L,\Gamma\rangle$, we have a fourth independent integral $F$ only in three remarkable cases (e.g., see \cite{BM}):

\begin{itemize}

\item\emph{Euler top}: $I=\diag(I_1,I_2,I_3)$, $\mathbf L=(L_1,L_2,L_3)$, $\chi=0$, $F=\langle \mathbf M+\mathbf L,\mathbf M+\mathbf L\rangle$, Zhykovskiy \cite{Zh} and  Volterra  \cite{Vol}.

\item\emph{Lagrange top}: $I=\diag(I_1,I_1,I_3)$, $\chi=(0,0,\chi_3)$, $\mathbf L=(0,0,\eta)$, $F=M_3$.

\item\emph{Kowalevski top}: $I=\diag(1,1,\frac12)$, $\chi=(\chi_1,0,0)$, $\mathbf L=(0,0,\eta)$,
\[
F=(M_1^2-M_2^2-2\chi_1\Gamma_1)^2+(2M_1M_2-2\chi_1\Gamma_2)^2+8\eta(M_3-2\eta)(M_1^2+M_2^2)-16\chi_1\eta M_1\Gamma_3.
\]

\end{itemize}

The Kowalevski integrable case (1889, \cite{Kow}) of a rigid body motion around a fixed point is one of the keystones of the classical
theory of integrable systems (e.g., see Dragovi\' c \cite{D2010} and reference therein).
Curiously, the Kowalevski top with the gyrostat is obtained 100 years later, almost simultaneously in three different papers by Yehia \cite{Y}, Komarov \cite{K}, and Reyman and  Semenov-Tian-Shansky \cite{RST}.

\section{Integrability of contact magnetic geodesic flows on $V_{n,2}$}\label{sec3}

 In the description of the magnetic Hamiltonian flows on $T^*V_{n,2}$, we follow the notation of \cite{FeJo_MZ_2012, FeJo_DCDS_2021}, where we have studied Neumann systems (see \cite{moser}) on Stiefel varieties.  The notation is adopted
for the rank two Stiefel variety and the phase space endowed with the twisted symplectic structure, the canonical one perturbed by the standard contact magnetic structure $\omega_{mag}$. For the Hamiltonian formalism of magnetic systems in a general setting, see
Novikov \cite{novikov}.

As usual, we consider vectors $x\in\R^n$ as columns and $x\wedge y=x\otimes y-y\otimes x=xy^T-yx^T\in so(n)$.
In what follows, we identify $so(n)\cong so(n)^*$  by the use of
the invariant metric $\langle \eta_1,\eta_2\rangle=-\frac12
\tr(\eta_1 \eta_2)$. The standard base in $\R^n$ is denoted by
\[
\mathbf{E}_1=(1,0,\dots,0,0), \, \mathbf E_2=(0,1,\dots,0,0),\, \dots, \, \mathbf{E}_n=(0,0,\dots,0,1).
\]

\subsection{Magnetic Hamiltonian systems on $T^*V_{n,2}$}\label{sec3a}

By the analogy to the constraints \eqref{e1e2}, \eqref{e-dot}, the cotangent bundle $T^*V(n,2)$
can be realized within $\R^{4n}(\e_1,\e_2,\p_1,\p_2)$ by the equations
\begin{equation}
\label{cxp}
\begin{aligned}
& f_{11}=\langle \e_1,\e_1\rangle =1, \quad f_{22}=\langle
\e_2,\e_2\rangle=1, \quad f_{12}=\langle \e_1,\e_2\rangle=0,\\
& g_{11}=\langle \e_1,\p_1\rangle =0, \quad g_{22}=\langle
\e_2,\p_2\rangle=0, \quad g_{12}=\langle \e_1,\p_2\rangle +\langle
\e_2,\p_1\rangle=0.
\end{aligned}
\end{equation}

The canonical symplectic structure $\omega$ on $T^*V_{n,2}$ is the
restriction of the canonical 2-form in the ambient space
$\R^{4n}(\e_1,\e_2,\p_1,\p_2)$:
\[
\omega=d\p_1\wedge d\e_1+d\p_2\wedge d\e_2\vert_{V_{n,2}}=\sum_{i=1}^n dp^i_1\wedge de^i_1+dp^i_2\wedge de^i_2\vert_{T^*V_{n,2}}.
\]

Let $\rho: T^*V_{n,2}\to V_{n,2}$ be the natural projection and let $\eta\in\R$ be a real parameter (strength
of the magnetic field). Define the twisted symplectic form  (see \cite{novikov})
\begin{equation}\label{twisted}
\omega_\eta=\omega+\eta\,\rho^*\omega_{mag}.
\end{equation}
By $\{\cdot,\cdot\}_\eta$ we denote the corresponding contact magnetic Poisson brackets on $T^*V_{n,2}$.

 Note that $(T^*V_{n,2},\omega_\eta)$ is a symplectic submanifold of $(\R^{4n},d\p_1\wedge d\e_1+d\p_2\wedge d\e_2+\eta\,d\e_1 \wedge d\e_2)$.
 Let $H(\e_1,\e_2,\p_1,\p_2)$ be a smooth function defined in $\R^{4n}$.
The Hamiltonian equations of the restriction $H\vert_{T^*V_{n,2}}$ can be obtained by using the Lagrange multipliers
as follows (e.g., see  pages 177 and 178 in Moser \cite{moser}
for the description of the Hamiltonian systems with constraints with respect to the standard symplectic form). Set
\[
H^*=H-\frac{\lambda_{11}}2 f_{11}-\lambda_{12}f_{12}-\frac{\lambda_{22}}2 f_{22}-\mu_{11}g_{11}-\mu_{12}g_{12}-\mu_{22}g_{22}.
\]

Consider the non-constrained Hamiltonian equations of $H^*$ on $\R^{4n}(\e_1,\e_2,\p_1,\p_2)$ with respect to the
twisted symplectic form
$d\p_1\wedge d\e_1+d\p_2\wedge d\e_2+\eta\,d\e_1 \wedge d\e_2$:
\begin{equation}\label{H*}
\begin{aligned}
\dot \e_1&=\frac{\partial H^*}{\partial \p_1}=\frac{\partial H}{\partial \p_1}-\mu_{11}\e_1-\mu_{12}\e_2,\\
\dot \e_2&=\frac{\partial H^*}{\partial \p_2}=\frac{\partial H}{\partial \p_2}-\mu_{12}\e_1-\mu_{22}\e_2,\\
\dot \p_1&=-\frac{\partial H^*}{\partial \e_1}+\eta\,\frac{\partial H^*}{\partial \p_2}\\
&=-\frac{\partial H}{\partial \e_1}+\eta\,\frac{\partial H}{\partial \p_2}
-\eta\big(\mu_{12}\e_1+\mu_{22}\e_2\big)+\lambda_{11}\e_1+\lambda_{12}\e_2+\mu_{11}\p_1+\mu_{12}\p_2,\\
\dot \p_2&=-\frac{\partial H^*}{\partial \e_2}-\eta\,\frac{\partial H^*}{\partial \p_1}\\
&=-\frac{\partial H}{\partial \e_2}-\eta\,\frac{\partial H}{\partial \p_1}+\eta\big(\mu_{11}\e_1+\mu_{12}\e_2\big)+\lambda_{12}\e_1+\lambda_{22}\e_2+\mu_{12}\p_1+\mu_{22}\p_2.
\end{aligned}
\end{equation}
The multipliers $\lambda_{ij}, \mu_{ij}$ are determined from the condition that the trajectories of the system satisfy the constraints \eqref{cxp}.
 In other words,
the Hamiltonian vector field of  $H\vert_{T^*V_{n,2}}$ on the symplectic submanifold $T^*V_{n,2}\subset\R^{4n}$,
at a point $(\e_1,\e_2,\p_1,\p_2)\in T^*V_{n,2}$, is the linear combination
 of the Hamiltonian vector fields of $H$ and the constraint functions ${f_{11}}$, $f_{12}$, $f_{22}$, ${g_{11}}$, $g_{12}$,
 $g_{22}$ in the ambient space $\R^{4n}$.

For a sake of simplicity of notation, when we consider Hamiltonian equations or Poisson brackets $\{\cdot,\cdot\}_\eta$ on
$(T^*V_{n,2},\omega_\eta)$ we will simply write $H$ instead of $H\vert_{T^*V_{n,2}}$ for the restriction of $H\in C^\infty(\R^{4n})$ to $T^*V_{n,2}$.

\subsection{$SO(n)$ and $SO(2)$--actions and submersion to the oriented Grassmmanian $G^+_{n,2}$}\label{sec3b}

Let $X=(\e_1,\e_2), P=(\p_1,\p_2)\in M_{n,2}(\R)$,  where $M_{n,2}(\R)$ denotes the space of real matrixes with $2$ columns and $n$ rows.
The Lie groups $SO(n)$ and $SO(2)$ naturally act on $T^*V(n,2)$ by left and right multiplications, respectively:
\begin{equation}\label{actions}
(X,P)\longmapsto (RX,RP), \quad (X,P)\longmapsto (XK^T,PK^T), \quad R\in SO(n), \quad K\in SO(2).
\end{equation}
Left $SO(n)$-action is a natural extension of the standard $SO(n)$--action on the Euclidean space $\R^n$, while right $SO(2)$-action defines rotations in the 2-plane
$\Span\{\e_1,\e_2\}$.  The following statement is proved in \cite{FeJo_MZ_2012}.

\begin{lem}
The left $SO(n)$ and right $SO(2)$ actions are Hamiltonian with the equivariant
magnetic momentum mappings given by
\begin{align}
& \Phi_\eta=PX^T-XP^T+\eta\, \e_1\wedge \e_2=\p_1\wedge \e_1+\p_2\wedge \e_2+\eta\, \e_1\wedge \e_2, \label{momentum_map_left} \\
& \Psi=\langle \e_1,\p_2\rangle -\langle \e_2,\p_1\rangle.
\label{moment_map_right}
\end{align}
\end{lem}

The Stiefel variety $V_{n,2}$ is diffeomorphic to the homogeneous space $SO(n)/SO(n-2)$ with respect to the left $SO(n)$--action. Namely, consider the point $o:=(\mathbf E_1,\mathbf E_2)\in V_{n,2}$. Then $V_{n,2}$ is the orbit of $o$: $SO(n)\cdot o=V_{n,r}$ and the isotropy group of $o$ is:
\[
SO(n)_o:=\{\diag(1,1,S)\,\vert\,  S\in SO(n-2)\}\cong SO(n-2).
\]
In particular, the contact distribution $\mathcal H$ is $SO(n)$--invariant:
\[
\mathcal H\vert_{(R\cdot o)}=R\cdot \mathcal H\vert_o, \qquad R\in SO(n).
\]

The Reeb vector field $Z$ on $(V_{n,2},\alpha)$, defined by the relations
$\alpha(Z)=1$, $i_Z\omega_{mag}=0$ (see \cite{LM}),
reads
\begin{equation}\label{reeb}
Z_{(\e_1,\e_2)}=(-\e_2,\e_1)\in T_{(\e_1,\e_2)} V_{(n,2)}.
\end{equation}

Therefore, the Reeb vector field induces right $SO(2)$--action on $V_{n,2}$ with its prolongation to $T^*V_{n,2}$ given by the Hamiltonian flow of the momentum mapping $\Psi$.
The quotient space $V_{n,2}/SO(2)$ is known as the Grassmmanian manifold $G^+_{n,2}$ of oriented $2$-planes in $\R^n$. The quotient mapping is given by
the submersion
\begin{equation}\label{submersion}
\pi\colon V_{n,2}\longrightarrow G^+_{n,2}, \qquad \pi(\e_1,\e_2)=\e_1\wedge \e_2.
\end{equation}

Thus, $(V_{n,2},\alpha)$ is an example of Boothby--Wang contact manifold \cite{BW}. The contact distribution $\mathcal H$ can be
seen as the horizontal space of the fibration \eqref{submersion}, and $\omega_{mag}$ is the basic 2-form proportional to the curvature of $\mathcal H$.
It induces Kirillov-Kostant-Souriau symplectic form on $G^+_{n,2}$ considered as the adjoint orbit of $\mathbf E_1\wedge\mathbf E_2$ in $so(n)$:
\[
G^+_{n,2} \cong \{ R\cdot \mathbf E_1\wedge \mathbf E_2 \cdot R^T\, \vert\, R\in SO(n)\} \cong SO(n)/SO(2)\times SO(n-2).
\]

\subsection{SO(n)--invariant magnetic geodesic flows}\label{sec3I}

The set of $SO(n)$--invariant functions on the cotangent bundles of Stiefel manifolds is described in Lemmas 2.3 and 2.4, \cite{FeJo_MZ_2012}.
For $n\ge 4$, there are 4  functionally independent functions in $C^\infty_{SO(n)}(T^*V_{n,2})$,  for example,
\begin{equation}\label{SO(n)}
\Psi, \qquad \langle \p_1,\p_1\rangle, \qquad \langle \p_2,\p_2\rangle, \qquad \langle \p_1,\p_2\rangle.
\end{equation}

Since $\langle \p_1,\e_2\rangle+\langle \p_2,\e_1\rangle=0\vert_{T^*V_{n,r}}$, instead of $\Psi$ we can take one of the polynomials
$\langle \p_1,\e_2\rangle$ or $\langle \p_2,\e_1\rangle$. Thus, for $n\ge 4$, a general form of the $SO(n)$--invariant Hamiltonian, quadratic in momenta, has the form
\begin{equation}
\label{quadratic}
H_{a}=\frac12 a_{1}\langle \p_1,\p_1\rangle+\frac12a_{2} \langle \p_2,\p_2\rangle+a_{3} \langle \p_1,\e_2\rangle\langle \p_2,\e_1\rangle+ a_{4} \langle \p_1,\p_2\rangle.
\end{equation}

For $n=3$ we have 3 linear in momenta independent functions on $T^*V_{3,2}$:
\begin{equation}\label{SO(3)}
\Psi, \qquad \langle \p_1, \e_1\times \e_2\rangle, \qquad \langle \p_2, \e_1\times \e_2\rangle,
\end{equation}
and a general form of the $SO(3)$--invariant Hamiltonian, quadratic in momenta, has the form
\begin{equation}
\label{quadratic3}
\begin{aligned}
H_{b}=&\frac12 b_{1}\langle \p_1,\e_1\times \e_2\rangle^2+\frac12 b_{2} \langle \p_2,\e_1\times \e_2\rangle^2+
\frac12 b_{3} \Psi^2 \\
&+b_4\langle \p_1,\e_1\times \e_2\rangle\langle \p_2,\e_1\times \e_2\rangle+b_{5}\Psi\langle \p_1,\e_1\times \e_2\rangle+ b_6\Psi\langle \p_2,\e_1\times \e_2\rangle.
\end{aligned}
\end{equation}

Assume that parameters $a=(a_1,a_2,a_3,a_4)$ ($b=(b_1,b_2,b_3,b_4,b_5,b_6)$ for $n=3$) determine $H_a$ (respectively, $H_b$) as a positive definite quadratic form in momenta on the cotangent bundle $T^*V_{n,2}$. By $ds^2_{a}$ ($ds^2_b$ for $n=3$) we denote the associate Riemannian metric on $V_{n,2}$.
Note that for $n=3$, the metrics $\{ds^2_a\}$ form a subclass in the set of all $SO(3)$-invariant metrics $\{ds^2_b\}$.

\begin{lem}\label{pozitivnost}
The parameters $a=(a_1,a_2,a_3,a_4)$ define $SO(n)$-invariant Riemannian metric $ds^2_a$ on $V_{n,2}$ if and only if
\begin{equation}\label{uslovi}
a_1>0, \quad a_2>0, \quad a_1+a_2>2a_3, \quad a_1a_2> a_4^2.
\end{equation}
\end{lem}

Among Riemannian metrics $\{ds^2_{a}\}$ we have two-parametric family of $SO(n)\times SO(2)$--invariant Riemannian metrics $\{ds^2_{\nu,\kappa}\}$ defined by
$(a_1,a_2,a_3,a_4)=(\nu,\nu,-\nu(1+2\kappa),0)$:
\begin{equation}\label{ham*}
H_{\nu,\kappa}=
 \frac{\nu}{2}\langle \p_1,\p_1\rangle +\frac{\nu}{2}\langle
\p_2,\p_2\rangle -\nu(1+2\kappa)\langle \p_1,\e_2\rangle \langle
\p_2,\e_1\rangle,
\end{equation}
where $\nu>0$ and $\kappa>-1$. In terms of the momentum mappings, we have:
\[
H_{\nu,\kappa}=\frac12\nu\langle \Phi_0,\Phi_0\rangle + \frac12 \nu\kappa\Psi^2.
\]

\begin{lem}\label{lezandr} The Lagrangian of the metric $ds^2_{\nu,\kappa}$, the Legandre transformation of the Hamiltonain $H_{\nu,\kappa}$ in the presence of the constraints \eqref{cxp} and \eqref{e1e2}, \eqref{e-dot}, is given by
\[
L_{\nu,\kappa}(\e_1,\e_2,\dot \e_1,\dot \e_2)=\frac1{2\nu}\langle \dot \e_1,\dot
\e_1\rangle+\frac1{2\nu}\langle \dot \e_2,\dot
\e_2\rangle+\frac{1+2\kappa}{2\nu(1+\kappa)}\langle \e_1,\dot
\e_2\rangle\langle \e_2,\dot \e_1\rangle.
\]
\end{lem}

For $\kappa=0$, $ds^2_{\nu,0}$ is a \emph{normal metric}, induced from a bi-invariant metric on the Lie group $SO(n)$ (it is unique up to multiplication by a positive constant). Further,  we have that the metric $ds^2_{1,-1/2}$ is the restriction to $V_{n,2}$ of the Euclidean metric in the ambient space $\R^{2n}(\e_1,\e_2)$.

\begin{lem}\label{rimanov-tok}
The Hamiltonian equations of the Hamiltonian \eqref{quadratic}
with respect to the twisted symplectic structure \eqref{twisted} are given by
\begin{equation}\label{geodezijski-tok}
\begin{aligned}
&\dot \e_1=a_1\p_1+a_4\p_2+a_3\langle \e_1,\p_2\rangle \e_2- a_4\langle \e_1,\p_2\rangle\e_1-\frac12\big(a_1\langle \e_2,\p_1\rangle+a_2\langle \e_1,\p_2\rangle\big)\e_2, \\
&\dot \e_2=a_2\p_2+a_4\p_1+a_3\langle \e_2,\p_1\rangle \e_1-\frac12\big(a_1\langle \e_2,\p_1\rangle+a_2\langle \e_1,\p_2\rangle\big) \e_1-a_4\langle \e_2,\p_1\rangle \e_2,\\
& \dot \p_1=-a_3\langle \p_1,\e_2\rangle \p_2+\eta\,\dot\e_2+ a_4\langle \e_1,\p_2\rangle\p_1+\frac12\big(a_1\langle \e_2,\p_1\rangle+a_2\langle \e_1,\p_2\rangle\big)\p_2+\lambda_{11} \e_1+ \lambda_{12} \e_2, \\
& \dot \p_2=-a_3\langle \p_2,\e_1\rangle \p_1-\eta\, \dot\e_1+\frac12\big(a_1\langle \e_2,\p_1\rangle+a_2\langle \e_1,\p_2\rangle\big) \p_1+a_4\langle \e_2,\p_1\rangle \p_2+\lambda_{12} \e_1+\lambda_{22} \e_2,
\end{aligned}
\end{equation}
where the Lagrange multipliers $\lambda_{ij}$ are
\begin{align*}
 & \lambda_{11}=({a_2-a_1})\langle \p_1,\e_2\rangle\langle \p_2,\e_1\rangle  -a_1 \langle \p_1,\p_1\rangle- a_4\langle \p_1,\p_2\rangle
 -\eta\,\big(\frac{a_1+a_2}{2}-a_3\big)\langle \p_2,\e_1\rangle, \\
 & \lambda_{12}=- 2a_4\langle \p_1,\e_2\rangle\langle \p_2,\e_1\rangle  -\frac{a_1+a_2}2 \langle \p_1,\p_2\rangle- \frac{a_4}2\langle \p_1,\p_1\rangle-\frac{a_4}2\langle \p_2,\p_2\rangle, \\
&\lambda_{22}=({a_1-a_2})\langle \p_1,\e_2\rangle\langle \p_2,\e_1\rangle  -a_2 \langle \p_2,\p_2\rangle- a_4\langle \p_1,\p_2\rangle
 +\eta\,\big(\frac{a_1+a_2}{2}-a_3\big)\langle \p_1,\e_2\rangle.
\end{align*}
\end{lem}

\noindent\emph{Proof.}
According to \eqref{H*}, the first two equations are
\begin{align*}
&\dot \e_1=a_1\p_1+a_3\langle \e_1,\p_2\rangle  \e_2+a_4\p_2-\mu_{11}\e_1-\mu_{12} \e_2, \\
& \dot \e_2=a_2\p_2+a_3\langle \e_2,\p_1\rangle \e_1+a_4\p_1-\mu_{12}\e_1-\mu_{22} \e_2.
\end{align*}
Since the solutions $(\e_1(t),\e_2(t),\p_1(t),\p_2(t))$ of the flow satisfy $f_{ij}(\e_1(t),\e_2(t))\equiv 0$, we get the equations
\begin{equation}\label{FM*}
\langle \e_1,\dot\e_1\rangle =0, \qquad \langle \e_2,\dot\e_2\rangle =0, \qquad \langle \e_1,\dot\e_2\rangle+\langle \dot \e_1,\e_2\rangle =0.
\end{equation}

By substituting the expressions for $\dot \e_1$ and $\dot\e_2$ into \eqref{FM*},
we get the multipliers $\mu_{ij}$,
\[
\mu_{11}=a_4\langle \e_1,\p_2\rangle, \qquad \mu_{22}=a_4\langle \e_2,\p_1\rangle, \qquad
\mu_{12}=\frac12\big(a_1\langle \e_2,\p_1\rangle+a_2\langle \e_1,\p_2\rangle\big),
\]
and the first two equations in \eqref{geodezijski-tok}.

Similarly, we have identities $g_{ij}(\e_1(t),\e_2(t),\p_1(t),\p_2(t))\equiv 0$ that imply
the equations
\begin{equation}\label{GL*}
\begin{aligned}
\langle \e_1,\dot\p_1\rangle+\langle \dot\e_1,\p_1\rangle =0, \qquad \langle \e_2,\dot\p_2\rangle+\langle \dot\e_2,\p_2\rangle =0, \\
\langle \e_1,\dot\p_2\rangle+\langle \dot\e_1,\p_2\rangle+\langle \e_2,\dot\p_1\rangle+\langle \dot\e_2,\p_1\rangle =0.
\end{aligned}
\end{equation}

Finally,
the Lagrange multipliers $\lambda_{ij}$ are determined from the equations obtained by substituting the last two equations in \eqref{geodezijski-tok}
and the derived expressions for $\dot\e_1$, $\dot\e_2$, into the relations \eqref{GL*}.
\hfill$\Box$

\begin{teo}\label{glavna}
The magnetic geodesic flows of the $SO(n)$--invariant metrics $\{ds^2_a\}$ (and $ds^2_b$ for $n=3$) are completely integrable in the non-commutative sense.
The complete algebra of first
integrals is generated by the components of the momentum mapping $\Phi_\eta^{ij}=\langle \Phi_\eta,\mathbf E_i\wedge \mathbf E_j\rangle$
and the Hamiltonian function.
For $n\ge 4$, the dimension of invariant isotropic manifolds is equal to $3$. For $n=3$, the dimension of invariant tori is equal to $2$.
\end{teo}

\noindent\emph{Proof.}
The proof essentially follows from the fact that $V_{n,2}=SO(n)/SO(n-2)$ is the complexity 1 homogeneous space, i.e., $(SO(n),SO(n-2))$ is an almost spherical pair (see \cite{MS}).
According to the Noether theorem, we have
\[
\{C^\infty_{SO(n)}(T^*V_{n,2}),\Phi_\eta^{ij}\}_\eta=0.
\]
Also, $SO(n)$--invariant functions on $T^*V_{n,2}$ together with the components of the momentum mapping $\Phi_\eta^{ij}$ generate a complete algebra of functions on $(T^*V_{n,2},\omega_\eta)$, see Theorem 1 in \cite{BJ2008}.

For $n\ge 4$, functionally independent functions in $C^\infty_{SO(n)}(T^*V_{n,2})$ are given by \eqref{SO(n)}.
Moreover, $SO(n)$-invariant functions
\[
J_1=\mathrm{tr}(\Phi_\eta^2), \qquad J_2=\mathrm{tr}(\Phi_\eta^4),
\]
are independent Casimirs within the algebra $(C^\infty_{SO(n)}(T^*V_{n,2}),\{\cdot,\cdot\}_\eta)$ and $\mathcal F=\{J_1, J_2, H_a\}$ is a maximal set of commuting functions within $(C^\infty_{SO(n)}(T^*V_{n,2}),\{\cdot,\cdot\}_\eta)$. Since $J_1$ and $J_2$ functionally depend on $\Phi_\eta^{ij}$, it follows that
\begin{equation}\label{integrali}
H_a, \qquad \Phi_\eta^{ij}, \qquad 1\le i<j\le n
\end{equation}
is a complete set of integrals on the symplectic manifold $(T^*V_{n,2},\omega_\eta)$.
Thus, the Hamiltonian flow is completely integrable in the non-commutative sense (see \cite{MF, N}). The regular invariant isotropic tori, level sets of the integrals
\eqref{integrali} are 3-dimensional. They are generated by the Hamiltonian flows of the functions $J_1, J_2, H_a$, which are in involutions with all
integrals:
\[
 \{H_a, \Phi_\eta^{ij}\}_\eta=0, \, \{J_1,J_2\}_\eta=0, \, \{J_k, \Phi_\eta^{ij}\}_\eta=0, \, \{H_a,J_k\}_\eta=0,\,  1\le i<j\le n, \, k=1,2.
\]

For $n=3$, $SO(3)$--invariant functions are given by \eqref{SO(3)}, we have only one independent Casimir function $J_1$ within
$(C^\infty_{SO(3)}(T^*V_{3,2}),\{\cdot,\cdot\}_\eta)$. The integrals
\begin{equation*}
H_b, \quad \Phi_\eta^{12}, \quad \Phi_{\eta}^{13}, \quad \Phi_\eta^{23}
\end{equation*}
are independent and their regular level sets are 2-dimensional isotropic tori. They are generated by the Hamiltonian flows of $H_b$ and $J_1$.
\hfill$\Box$

\begin{rem}[Liouville integrability]
The magnetic geodesic flows of the metrics $\{ds^2_a\}$ are completely integrable in the Liouville sense as well. For example, for $n\ge 4$, by using the chain of
subalgebras $so(2)<so(3)<\dots<so(n)$ (see \cite{Thimm, JSV2023}),  we get $2n-3$ independent commuting integrals
\begin{align*}
&F_1=H_a, \quad F_2=\Phi_{\eta}^{12}, \quad  F_{k}=\tr\big(\sum_{1\le i<j\le k} \Phi_\eta^{ij} \mathbf E_i\wedge \mathbf E_j\big)^2, \quad k=3,\dots,n, \\
& F_{n+l}=\tr\big(\sum_{1\le i<j\le 3+l} \Phi_\eta^{ij} \mathbf E_i\wedge \mathbf E_j\big)^4, \quad l=1,\dots,n-3,  \quad (\text{note that $F_n=J_1$, $F_{2n-3}=J_2$}).
\end{align*}
For $n=3$, commuting integrals are $F_1=H_b$, $F_2=\Phi_{\eta}^{12}$, and $F_3=J_1$.
\end{rem}

For $\eta=0$, Theorem \ref{glavna} is a generalization of Theorem 4.5 \cite{FeJo_MZ_2012}, where it is proved that $SO(n)\times SO(r)$--invariant metrics
on rank $r$ Stiefel manifolds $V_{n,r}$ have integrable geodesic flows. For $r>2$, the complexity of the homogeneous spaces $V_{n,r}=SO(n)/SO(n-r)$ is greater than 1, and the proof of Theorem \ref{glavna} can not be extended to a generic $SO(n)$--invariant metric on $V_{n,r}$.

Other examples of integrable magnetic geodesic flows on Lie groups and homogeneous spaces can be found in \cite{AS2020, BJ2006, BJ2008, E2005, MSY2008}.

\begin{rem}[Manakov metrics]\label{primedba}
We have integrability of the geodesic flows for
the class of the singular Manakov metrics on $SO(n)$ and the Manakov metrics on homogeneous spaces of the Lie group $SO(n)$ (see \cite{DGJ2009}), in particular
the rank $r$ Stiefel varieties $V_{n,r}=SO(n)/SO(n-r)$, established in \cite{DGJ2009, DGJ2015, Myk}.
The Manakov metrics on the rank 2 Stiefel variety $V_{n,2}$ are defined as follows. The tangent space at $o=(\mathbf E_1,\mathbf E_2)$ can be naturally identified with the subspace of the Lie algebra $so(n)$ orthogonal to the Lie algebra $so(n-2)$ of the isotropy group $SO(n)_o\cong SO(n-2)$:
$T_oV_{n,2}\cong \mathfrak v=\Span\{\mathbf E_i\wedge\mathbf E_j\,\vert\, i<j, \, 1\le i \le 2, \, 2\le j \le n \}$.
We have the decomposition
\begin{equation}\label{dekompozicija}
T_o V_{n,2}\cong \mathfrak v=\mathfrak v_{1,2}\oplus \mathfrak v_{1,3}\oplus \mathfrak v_{2,3},
\end{equation}
where
\[
\mathfrak v_{1,2}=\Span\{\mathbf E_1\wedge\mathbf E_2\}, \, \mathfrak v_{1,3}=\Span\{\mathbf E_1\wedge\mathbf E_j\,\vert\, 3\le j \le n \},\,
\mathfrak v_{2,3}=\Span\{\mathbf E_2\wedge\mathbf E_j\,\vert\, 3\le j \le n \}.
\]
Now, for diagonal $n\times n$ matrixes
$A=(\alpha_1,\alpha_2,\alpha_3,\dots,\alpha_3)$, $B=(\beta_1,\beta_2,\beta_3,\dots,\beta_3)$,
the Manakov metric $ds^2_{A,B}$ is defined as the $SO(n)$--invariant metric on $V_{n,2}$ with the restriction to $\mathfrak v$ given by
(see eq. (40) in \cite{DGJ2009}):
\begin{equation}\label{manakov}
\langle \xi,\eta\rangle_{A,B}=\frac{\alpha_1-\alpha_2}{\beta_1-\beta_2}\langle \xi_{1,2},\eta_{1,2}\rangle
+\frac{\alpha_1-\alpha_3}{\beta_1-\beta_3}\langle \xi_{1,3},\eta_{1,3}\rangle+
\frac{\alpha_2-\alpha_3}{\beta_2-\beta_3}\langle \xi_{2,3},\eta_{2,3}\rangle,
\end{equation}
where $\xi=\xi_{1,2}+\xi_{1,3}+\xi_{2,3}$, $\eta=\eta_{1,2}+\eta_{1,3}+\eta_{2,3}$, $\xi_{i,j}, \eta_{i,j}\in\mathfrak v_{i,j}$.
It can be proved that the parameters $a=(a_1,a_2,a_3,a_4)$ in the Hamiltonian \eqref{quadratic} and the matrixes $A, B$ are related by
\begin{equation}\label{parametri-manakov}
a_1=\frac{\beta_1-\beta_3}{\alpha_1-\alpha_3}, \,\, a_2=\frac{\beta_2-\beta_3}{\alpha_2-\alpha_3},
\,\, a_3=\frac12\big(\frac{\beta_1-\beta_3}{\alpha_1-\alpha_3}+\frac{\beta_2-\beta_3}{\alpha_2-\alpha_3}\big)-2\frac{\beta_1-\beta_2}{\alpha_1-\alpha_2}, \,\, a_4=0.
\end{equation}
Furthermore,  the conditions \eqref{uslovi} are equivalent to the obvious conditions
\[
\frac{\beta_1-\beta_3}{\alpha_1-\alpha_3}>0, \qquad \frac{\beta_2-\beta_3}{\alpha_2-\alpha_3}>0, \qquad \frac{\beta_1-\beta_2}{\alpha_1-\alpha_2}>0
\]
on the parameters $\alpha_i,\beta_i$ so that \eqref{manakov} is a positive definite scalar product.
\end{rem}

\subsection{$n=3$ and a free rigid body with a gyrostat}\label{sec3III}

The Stiefel variety $V_{3,2}$ is naturally isomorphic to $SO(3)$:
\[
R: V_{3,2}\to SO(3): \quad R(\e_1,\e_2)=(\e_1\, \e_2\, \e_3), \quad
\e_3=\e_1\times \e_2.
\]

Within this identification, we can see $V_{3,2}$ as the
configuration space of a rigid body motion around a fixed point.
The matrix $R$ maps the frame attached to the body to
the fixed reference frame:
\[
\e_1=R(\e_1,\e_2)\mathbf E_1, \quad \e_2=R(\e_1,\e_2) \mathbf E_2, \quad
\e_3=R(\e_1,\e_2)\mathbf E_3.
\]

It is well known  that the momentum mapping $\Phi_0$ of the left $SO(3)$-action on $T^*SO(3)$ corresponds to the momentum $\mathbf m=(m_1,m_2,m_3)$
of the rigid body in the fixed reference frame,
\[
\Phi_0=\p_1\wedge \e_1+\p_2\wedge \e_2=
\begin{pmatrix}
0 & -m_3 & m_2 \\
m_3 & 0 & -m_1 \\
-m_2 & m_1 & 0
\end{pmatrix}
\]
while the matrix
\[
R^T \Phi_0 R=
\begin{pmatrix}
0 & \langle \p_2,\e_1\rangle - \langle \p_1,\e_2\rangle & -\langle
\p_1,\e_3\rangle \\
 \langle \p_1,\e_2\rangle-\langle \p_2,\e_1\rangle  & 0 & -\langle
\p_2,\e_3\rangle \\
\langle \p_1,\e_3\rangle  & \langle \p_2,\e_3\rangle & 0
\end{pmatrix}
=
\begin{pmatrix}
0 & -M_3 & M_2 \\
M_3 & 0 & -M_1 \\
-M_2 & M_1 & 0
\end{pmatrix}
\]
corresponds to the momentum $\mathbf M=(M_1,M_2,M_3)$ of a rigid body in the moving reference
frame (we follows Arnold's notation \cite{Ar}).

Thus,  we can write the Hamiltonian $H_b$, and in particular $H_{\nu,\kappa}$, in the usual left-invariant form on $SO(3)$ as
\begin{equation}\label{general-rigid-body}
H_{b}=\frac12 b_{1} M_2^2+\frac12 b_{2} M_1^2+
\frac12 b_{3} M_3^2 -b_4 M_1M_2 +b_{5} M_2M_3- b_6 M_1M_3
\end{equation}
and
\begin{equation}\label{rigid-body}
H_{\kappa,\nu}=\frac12\nu(M_1^2+M_2^2+(1+\kappa)M_3^2).
\end{equation}

The Hamiltonian functions \eqref{general-rigid-body} and \eqref{rigid-body} model the motion of a free rigid body around a fixed point with inertia operators given respectively by
\begin{equation}\label{IO-gen}
I_b=A_b^{-1}, \qquad
A_b=\begin{pmatrix}
b_2 & -b_4 & -b_6 \\
-b_4 & b_1 &  b_5 \\
-b_6 & b_5 & b_3
\end{pmatrix},
\end{equation}
and
\begin{equation}\label{IO}
I_{\nu,\kappa}=\nu^{-1}\diag\big(1,1,\frac{1}{1+\kappa}\big).
\end{equation}

We used the above identification to perform discretization of the Euler top with the Euclidean metric ($\nu=1$, $\kappa=-1/2$) in \cite{FeJo_PSIM_2020}.

Similarly, the matrix
\[
R^T (\eta \e_1\wedge \e_2) R=\eta E_1\wedge E_2
=
\begin{pmatrix}
0 & \eta & 0 \\
-\eta & 0 & 0 \\
0    & 0  & 0
\end{pmatrix}
\]
corresponds to the constant vector $\mathbf L=(0,0,-\eta)$ in the moving reference frame.
Let $\mathbf l$ be the vector $\mathbf L$ considered in the fixed frame.
The conservation of the momentum $\Phi_\eta$ (i.e., conservation of $\mathbf m+\mathbf l$) corresponds to the Euler equations in the moving reference frame (see \cite{Ar})
\begin{equation}\label{VZ}
\frac{d}{dt}\big(\mathbf M+\mathbf L\big)=(\mathbf M+\mathbf L)\times \Omega, \qquad \Omega=A_b\mathbf M.
\end{equation}

The system \eqref{VZ} represents the equations of the body with a gyrostat \eqref{MEP} without the potential forces.
We can summarize the above considerations in the following statement.

\begin{prp}\label{stav1}
For $n=3$, the magnetic geodesic flows of the metric $ds^2_b$ represents inertial motion of a system rigid body+gyrostat  around a fixed point.  The inertia operator
of the body is given by \eqref{IO-gen} and the gyrostat momentum is $\mathbf L=(0,0,-\eta)$
\end{prp}

It is clear that by taking the base of principal axes of inertia, we get the general form of the Zhukovskiy-Volterra gyrostat with $I=\diag(I_1,I_2,I_3)$ and
$\mathbf L=(L_1,L_2,L_3)$.  Further, the magnetic geodesic flow of the metric $ds^2_{\nu,\kappa}$ represents the $SO(2)$--symmetric Zhukovskiy-Volterra gyrostat.

 The Hamiltonian \eqref{rigid-body}, for $\kappa=-1$, defines
the standard  sub-Riemannian geodesic structure on $SO(3)$ (see Section \ref{secSR}).

\begin{rem}
The sphere $S^3 \cong SU(2)$ is a two-covering of the group $SO(3)\cong V_{3,2}$. The contact
form $\alpha$ and the metric $ds^2_{\nu,\kappa}$ on $V_{3,2}$ correspond to the standard contact form and the $S^3\times SO(2)$--invariant Riemannian metric on $S^3$, which we also denote by $\alpha$ and $ds^2_{\nu,\kappa}$.
The sphere $(S^3,ds^2_{\nu,\kappa})$ is known as the Berger sphere.
The geometry of the contact magnetic geodesics on the Berger sphere has recently been described in details \cite{IM2025}.
Magnetic geodesic flows defined by exact homogeneous magnetic fields on round spheres of arbitrary dimensions are studied in
\cite{DGJ2023, DGJ2025, BKM2025, DGJ2025b}. They are completely integrable and closely related to the Neumann systems  \cite{BKM2025, DGJ2025b}.
Note that the standard contact form on
$S^{2m+1}=\{\mathbf x\in \R^{2m+2}\,\vert \, \langle \mathbf x,\mathbf x\rangle=1\}$, in particular on $S^3$, is defined by $\alpha=x_1dx_2+\dots+x_{2m+1}dx_{2m+2}$.
As in Theorem \ref{glavna} applied to $V_{3,2}$, the magnetic geodesic flows on $S^{2m+1}$ with respect to the contact magnetic fields
$\eta\cdot d\alpha$ ($\eta\ne 0$) are completely integrable in the non-commutative sense with two-dimensional generic invariant isotropic tori (see Theorem 5 in \cite{DGJ2025}).
\end{rem}

\section{Sub-Riemannian magnetic geodesic flows on $V_{n,2}$}\label{secSR}

Note that the decomposition \eqref{dekompozicija} is invariant with respect to
the adjoint $SO(n)_o\cong SO(n-2)$--action on $\mathfrak v$.
Thus, from the decomposition, we obtain the following $SO(n)$--invariant distributions on $V_{n,2}=SO(n)/SO(n-2)$, with
the restriction to $\mathfrak v=T_o V_{n,2}$ given by
\[
\mathfrak h=\mathfrak v_{1,3}\oplus \mathfrak v_{2,3}, \quad \mathfrak d_0=\mathfrak v_{1,2}\oplus \mathfrak v_{1,3},
\quad {\mathfrak d_{\frac{\pi}2}}=\mathfrak v_{1,2}\oplus \mathfrak v_{2,3}, \quad \mathfrak g_0=\mathfrak v_{1,3}, \quad \mathfrak g_{\frac{\pi}2}=\mathfrak v_{2,3}.
\]

The first one is the contact distribution $\mathcal H\vert_o=\mathfrak h$ (note that $\mathfrak v_{1,2}$ is the direction of the Reeb vector field $Z\vert_o$).
By $\mathcal D_{0}$, $\mathcal D_{\frac{\pi}2}$, $\mathcal G_0$, $\mathcal G_{\frac{\pi}2}$ we denote, the invariant distributions corresponding, respectively, to
$\mathfrak d_0$, $\mathfrak d_{\frac{\pi}2}$, $\mathfrak g_0$, $\mathfrak g_{\frac{\pi}2}$.

In the standard representation of the tangent bundle \eqref{e1e2}, \eqref{e-dot}, we have:
\begin{align*}
&\mathcal H\vert_{(\e_1,\e_2)}=(\Span\{\e_1,\e_2\}^\perp,\Span\{\e_1,\e_2\}^\perp),\\
&\mathcal D_0\vert_{(\e_1,\e_2)}=\Span\{Z\vert_{(\e_1,\e_2)}, \mathcal G_0\vert_{(\e_1,\e_2)}\},\\
&\mathcal D_{\frac{\pi}2}\vert_{(\e_1,\e_2)}=\Span\{Z\vert_{(\e_1,\e_2)}, \mathcal G_{\frac\pi2}\vert_{(\e_1,\e_2)}\},\\
&\mathcal G_0\vert_{(\e_1,\e_2)}=\{(\xi-\langle \xi,\e_2\rangle\e_2,0)\,\vert\, \xi\in\R^n, \langle \xi,\e_1\rangle=0\},\\
&\mathcal G_{\frac{\pi}2}\vert_{(\e_1,\e_2)}=\{(0,\xi-\langle \xi,\e_1\rangle\e_1)\,\vert\, \xi\in\R^n, \langle \xi,\e_2\rangle=0\}.
\end{align*}

The distributions $\mathcal H$, $\mathcal D_0$ and $\mathcal D_{\frac{\pi}2}$ are bracket generating, while $\mathcal G_0$ and $\mathcal G_{\frac{\pi}2}$
are not bracket generating distributions. Since the right $SO(2)$--action commute with the left $SO(n)$--action, other invariant bracket generating distributions
can be obtained using the rotations
\[
\sigma_\theta\colon (\e_1,\e_2) \longmapsto (\cos\theta\e_1+\sin\theta\e_2,-\sin\theta \e_1+\cos\theta \e_2),
\]
as
$
\mathcal D_\theta =d\sigma_\theta (\mathcal D_0).
$
Specifically, $\mathcal D_{\frac\pi2}$ is the image of $\mathcal D_0$ with respect to $\sigma_{\frac\pi2}$ (the involution in the projectivization of the two-plane $\Span\{\e_1,\e_2\}$).

Therefore, $\mathcal D_\theta\cong \mathcal D_0$ and
it is sufficient to describe sub-Riemannian structures on $\mathcal D_0$.

For $n=3$, we have the additional isomorphism $\mathcal H\cong\mathcal D_\theta$.

%The second and the third distributions are equivalent
%($\mathfrak d$ and
%$\mathfrak d$ are the orthogonal complements with respect to the killing form of two subalgebras in $so(n)$ isomorphic to $so(n-1)$).
%We consider below the distribution $\mathcal D$ defined by
%\[
%\mathcal D\vert_o=\mathfrak v_{1,2}\oplus \mathfrak v_{1,3}.
%\]

 All the statements obtained in this section, in the limit $\eta=0$, are valid for sub-Riemannian geodesics
without the influence of the magnetic field.

\subsection{Sub-Riemannian structures on the contact distribution $\mathcal H$}\label{secSRH}

According to Lemma \ref{lezandr}, the metric $ds^2_{\nu,\kappa}$ is defined by
\[
\langle \xi,\eta\rangle_{\nu,\kappa}=\frac{1}{\nu}\langle\xi_1,\eta_1\rangle+\frac{1}{\nu}\langle\xi_2,\eta_2\rangle+
 \frac{1+2\kappa}{2\nu(1+\kappa)}\big(\langle \e_1,\eta_2\rangle\langle \e_2,\xi_1\rangle+\langle \e_1,\xi_2\rangle\langle \e_2,\eta_1\rangle\big),
\]
$\xi=(\xi_1,\xi_2), \eta=(\eta_1,\eta_2)\in T_{(\e_1,\e_2)}V_{n,2}$. Therefore, the Reeb vector field \eqref{reeb} is orthogonal to the contact distribution $\mathcal H$ with respect to  $ds^2_{\nu,\kappa}$.
Moreover, we have
\[
\| Z \|_{\nu,\kappa}=
\sqrt{\frac{1}{\nu(1+\kappa)}}.
\]

The metric $ds^2_{\nu,\kappa}$ becomes singular as $\kappa$ tends to $-1$. However, if we consider the scalar product
$\langle \cdot,\cdot \rangle_{\nu,\kappa}$ restricted to  the contact distribution $\mathcal H$ it does not depend on $\kappa$:
\begin{equation}\label{restrikcija}
\langle \xi,\eta\rangle_{\nu,\kappa}=\frac{1}{\nu}\langle\xi_1,\eta_1\rangle+\frac{1}{\nu}\langle\xi_2,\eta_2\rangle, \qquad \xi,\eta\in\mathcal H\vert_{(\e_1,\e_2)}.
\end{equation}

\begin{figure}[h]
{\centering
{\includegraphics[width=10cm]{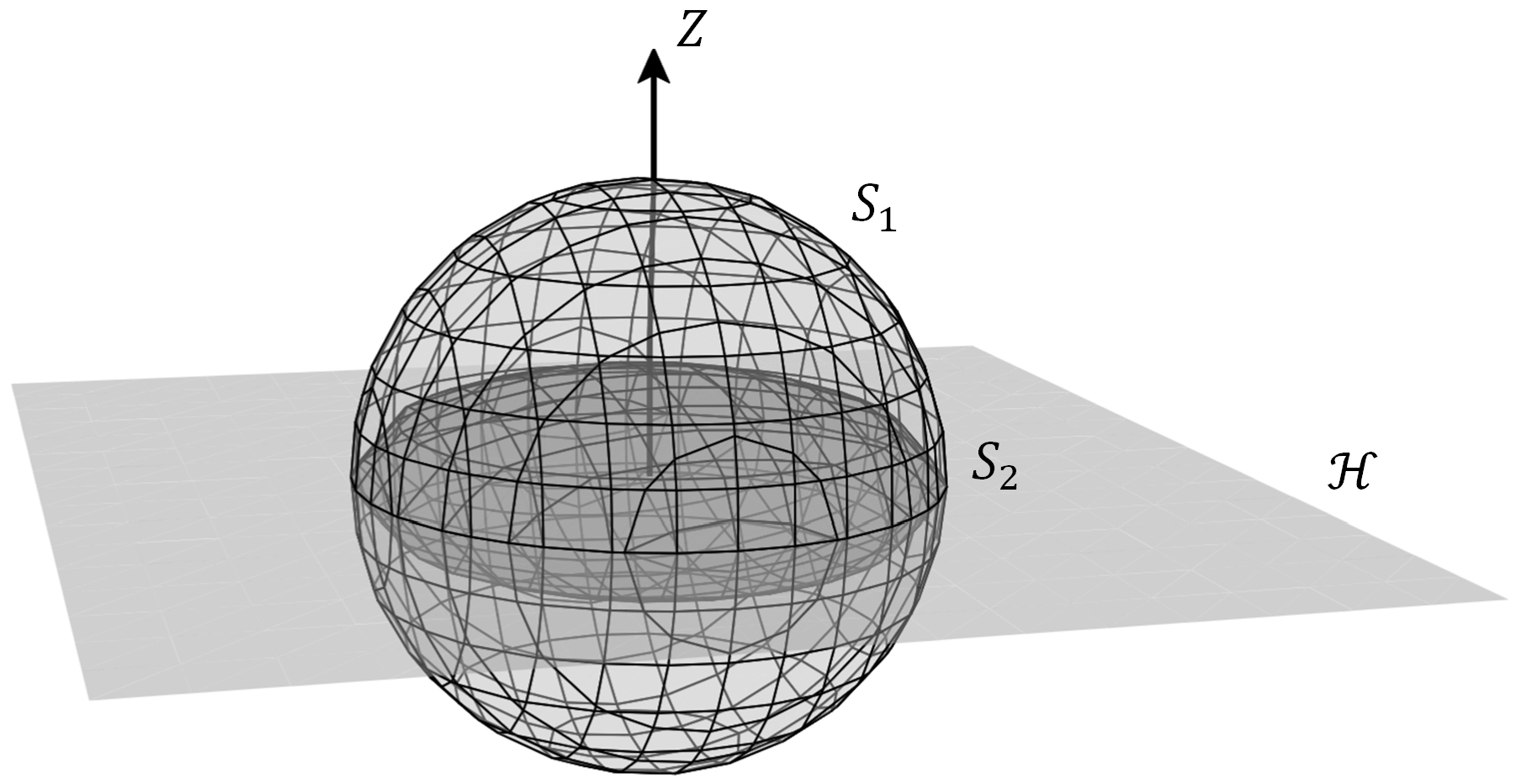}}
\caption{The unit spheres $S_1=\{\|\xi\|_{1,-\frac{1}{2}}=1\}$ and $S_2=\{\|\xi\|_{1,-\frac{15}{16}}=1\}$
within the tangent space $T_{(\mathbf e_1,\mathbf e_2)}V_{n,2}=\mathcal H\vert_{(\e_1,\e_2)}+\Span \{Z\vert_{(\e_1,\e_2)}\}$ with respect to the metrics $ds^2_{1,-\frac{1}{2}}$ and $ds^2_{1,-\frac{15}{16}}$.
Note that $\|Z\|_{1,-\frac{1}{2}}=\sqrt{2}$ and $\|Z\|_{1,-\frac{15}{16}}=4$.} \label{SRM}
}
\end{figure}

 Let $ds^2_{\mathcal H,\nu}$ be the sub-Riemannian structure on $\mathcal H$ defined by
\eqref{restrikcija}. In particular, $ds^2_{\mathcal H}:=ds^2_{\mathcal H,1}$
is the restriction of the Euclidean metric from
$\R^{2n}(\e_1,\e_2)$ to $\mathcal H$   and we refer to it as
the \emph{standard sub-Riemannian structure} on $V_{n,r}$. Note that the unit sphere bundle $\{\|\xi\|_{1,\kappa}\}=1\subset TV_{n,2}$, tends to the unit sphere bundle
within $\mathcal H$, as $\kappa$ tends to $-1$ (see Fig. \ref{SRM}).
On the other hand, the limit $\kappa\to -1$ is well defined for the Hamiltonian function $H_{1,\kappa}$, which represents
the Hamiltonian of the geodesic flow of the standard sub-Riemannian metric $ds^2_{\mathcal H}$:
\begin{equation*}
H_{\mathcal H}=H_{1,-1}=\frac12\langle \p_1,\p_1\rangle +\frac12\langle
\p_2,\p_2\rangle +\langle \p_1,\e_2\rangle \langle
\p_2,\e_1\rangle. \label{HamSR}
\end{equation*}

 The Riemannian metrics $ds^2_{1,\kappa}$ tame  $ds^2_{\mathcal H}$ and the above limit is usual
in construction of Hamiltonian functions of sub-Riemannian geodesic flows
from Hamiltonian functions of associated Riemannian metrics (e.g, see \cite{M2002}, page 32 and \cite{bloch}, page 389).

Similarly, in the limit $a_1+a_2 \to 2a_3$, for the parameters
\[
a_1>0, \quad a_2>0, \quad a_1+a_2=2a_3, \quad a_1a_2> a_4^2,
\]
we get the Hamiltonian
\begin{equation}\label{HamSRa}
H_{\mathcal H,a}=\frac12 a_{1}\langle \p_1,\p_1\rangle+\frac12a_{2} \langle \p_2,\p_2\rangle+\frac12(a_1+a_2)\langle \p_1,\e_2\rangle\langle \p_2,\e_1\rangle+ a_{4} \langle \p_1,\p_2\rangle,
\end{equation}
which defines the sub-Riemannian structure on $\mathcal H$ that we shall denote by $ds^2_{\mathcal H,a}$. For $a_1=1$, $a_2=1$, $a_4=0$ we get the standard
sub-Riemannian metric $ds^2_{\mathcal H}$.

From Lemma \ref{rimanov-tok}, for $a_1+a_2=2a_3$, we get.

\begin{lem}\label{sub-rimanov-tok}
(i) The Hamiltonian equations of the Hamiltonian \eqref{HamSRa}
with respect to the twisted symplectic structure \eqref{twisted} are given by:
\begin{equation}\label{MSRF}
\begin{aligned}
&\dot \e_1=a_1\p_1+a_4\p_2-a_4\langle \e_1,\p_2\rangle \e_1+a_1\langle \e_1,\p_2\rangle \e_2, \\
&\dot \e_2=a_2\p_2+a_4\p_1-a_4\langle \e_2,\p_1\rangle \e_2+a_2\langle \e_2,\p_1\rangle \e_1,\\
& \dot \p_1=-a_2\langle \p_1,\e_2\rangle \p_2+\eta\,\dot\e_2+ a_4\langle \e_1,\p_2\rangle\p_1+\lambda_{11} \e_1+ \lambda_{12} \e_2, \\
& \dot \p_2=-a_1\langle \p_2,\e_1\rangle \p_1-\eta\, \dot\e_1+a_4\langle \e_2,\p_1\rangle \p_2+\lambda_{12} \e_1+\lambda_{22} \e_2,
\end{aligned}
\end{equation}
where the Lagrange multipliers $\lambda_{ij}$ are
\begin{align*}
 & \lambda_{11}=({a_2-a_1})\langle \p_1,\e_2\rangle\langle \p_2,\e_1\rangle  -a_1 \langle \p_1,\p_1\rangle- a_4\langle \p_1,\p_2\rangle, \\
 & \lambda_{12}=- 2a_4\langle \p_1,\e_2\rangle\langle \p_2,\e_1\rangle  -\frac{a_1+a_2}2 \langle \p_1,\p_2\rangle- \frac{a_4}2\langle \p_1,\p_1\rangle-\frac{a_4}2\langle \p_2,\p_2\rangle, \\
&\lambda_{22}=({a_1-a_2})\langle \p_1,\e_2\rangle\langle \p_2,\e_1\rangle  -a_2 \langle \p_2,\p_2\rangle- a_4\langle \p_1,\p_2\rangle.
\end{align*}

(ii) The magnetic sub-Riemannian geodesic flows of the standard metrics $ds^2_{\mathcal H}$ is given by the equations:
\begin{equation}\label{MSF}
\begin{aligned}
&\dot \e_1=\p_1+\langle \e_1,\p_2\rangle  \e_2, \\
& \dot \e_2=\p_2+\langle \e_2,\p_1\rangle \e_1, \\
& \dot \p_1=-\langle \p_1,\e_2\rangle \p_2+\eta\,(\p_2+\langle \e_2,\p_1\rangle \e_1)-\langle \p_1,\p_1\rangle \e_1-\langle \p_1,\p_2\rangle \e_2, \\
& \dot \p_2=-\langle \p_2,\e_1\rangle \p_1-\eta\, (\p_1+\langle \e_1,\p_2\rangle  \e_2)-\langle \p_1,\p_2\rangle \e_1-\langle \p_2,\p_2\rangle \e_2.
\end{aligned}
\end{equation}
\end{lem}

Note  that $(\dot \e_1,\dot \e_2)\in\mathcal H_{(\e_1,\e_2)}=\big(\Span\{\e_1,\e_2\}^\perp,\Span\{\e_1,\e_2\}^\perp\big)$, and
the projections
$\rho(\e_1(t),\e_2(t),\p_1(t),\p_2(t))=(\e_1(t),\e_2(t))$ of the solutions of the flow \eqref{MSRF}, as required, are allowed by the contact distribution.

\begin{rem}
From the Pontryagin Maximum Principle (see \cite{Ju1996, AS}), it follows that apart from normal sub-Riemannian geodesics described by Hamiltonian equations there may exist abnormal geodesics, which depend only on the distribution.
Despite the existence and study of abnormal geodesics in one of the central problems in sub-Riemannian geometry  (e.g., see \cite{ABB, M2002, PS2}), here we focus on normal geodesics.  In the case of a contact distribution there are no abnormal
geodesics (see \cite{ABB} for the case of 3-dimensional contact manifolds). Therefore, the geodesics described in Lemma \ref{sub-rimanov-tok}
exhaust all geodesics for the corresponding sub-Riemannian
structure.
\end{rem}

\subsection{Sub-Riemannian structures on $\mathcal D_0$}\label{secSRD}

We can consider the limit $a_2, a_4 \to 0$
in the domain of parameters \eqref{uslovi},
to obtain the sub-Riemannian structures $ds^2_{\mathcal D_0,a}$ on the distribution $\mathcal D_0$ associated to the Hamiltonian function
\begin{align*}
H_{\mathcal D_0,a}=\frac12 a_{1}\langle \p_1,\p_1\rangle+a_{3} \langle \p_1,\e_2\rangle\langle \p_2,\e_1\rangle,  \qquad a_1>0,  \qquad a_1>2a_3.
\end{align*}

The restriction of the Euclidean metric from the ambient space $\R^{2n}(\e_1,\e_2)$ to the distribution $\mathcal D_0$ defines the \emph{standard sub-Riemannian structure} $ds^2_{\mathcal D_0}$.

\begin{lem}
The Hamiltonian function of the standard sub-Riemannian structure $ds^2_{\mathcal D_0}$ is:
\[
H_{\mathcal D_0}=\frac12 \langle \p_1,\p_1\rangle-\frac12\langle \p_1,\e_2\rangle\langle \p_2,\e_1\rangle.
\]
\end{lem}

\begin{lem}\label{sub-rimanov-tok*}
(i) The equations of the magnetic sub-Riemannian geodesic flow of the metric $ds^2_{\mathcal D_0,a}$
with respect to the twisted symplectic structure \eqref{twisted} are given by:
\begin{equation*}
\begin{aligned}
&\dot \e_1=a_1\p_1+\big(a_3+\frac{a_1}2\big)\langle \e_1,\p_2\rangle \e_2, \\
&\dot \e_2=\big(a_3-\frac{a_1}2\big)\langle \e_2,\p_1\rangle \e_1,\\
& \dot \p_1=-\big(a_3-\frac{a_1}2\big)\langle \p_1,\e_2\rangle \p_2+\eta\,\big(a_3-\frac{a_1}2\big)\langle \e_2,\p_1\rangle \e_1+\lambda_{11} \e_1+ \lambda_{12} \e_2, \\
& \dot \p_2=-\big(a_3+\frac{a_1}2\big)\langle \p_2,\e_1\rangle \p_1-\eta\, \big(a_1\p_1+\big(a_3+\frac{a_1}2\big)\langle \e_1,\p_2\rangle \e_2\big)+\lambda_{12} \e_1+
\lambda_{22} \e_2,
\end{aligned}
\end{equation*}
where the Lagrange multipliers $\lambda_{ij}$ are
\begin{align*}
 & \lambda_{11}=-a_1 \langle \p_1,\e_2\rangle\langle \p_2,\e_1\rangle  -a_1 \langle \p_1,\p_1\rangle-\eta\, \big(a_3-\frac{a_1}2\big)\langle \p_1,\e_2\rangle, \\
 & \lambda_{12}= -\frac{a_1}2 \langle \p_1,\p_2\rangle, \\
&\lambda_{22}=a_1\langle \p_1,\e_2\rangle\langle \p_2,\e_1\rangle +\eta\,\big(a_3-\frac{a_1}2\big) \langle \p_2,\e_1\rangle.
\end{align*}

(ii) Specially, the magnetic geodesic flow of the standard sub-Riemannian structure $ds^2_{\mathcal D_0}$ is given by
the equations:
\begin{equation*}
\begin{aligned}
&\dot \e_1=\p_1, \\
&\dot \e_2=-\langle \e_2,\p_1\rangle \e_1,\\
& \dot \p_1=\langle \p_1,\e_2\rangle \p_2-\big(\langle \p_1,\e_2\rangle\langle \p_2,\e_1\rangle+\langle \p_1,\p_1\rangle\big) \e_1-\frac{1}2 \langle \p_1,\p_2\rangle\e_2, \\
& \dot \p_2=-\eta\, a_1\p_1-\frac{1}2\langle \p_1,\p_2\rangle \e_1+
\big(\langle \p_1,\e_2\rangle\langle \p_2,\e_1\rangle -\eta\, \langle \p_2,\e_1\rangle \big)\e_2.
\end{aligned}
\end{equation*}
\end{lem}

\subsection{Integrability of magnetic sub-Riemannian geodesic flows}\label{sec3II}

Denote the considered sub-Riemannian structures
by $ds^2_{sR,a}$ and the corresponding Hamiltonian function by $H_{sR,a}$. For $n=3$ we only consider $ds^2_{\mathcal H,a}$.

With the same proof as Theorem \ref{glavna},  where instead of the Hamiltonian functions of $SO(n)$--invariant Riemannian metrics we consider
the Hamiltonian functions of $SO(n)$--invariant sub-Riemannian metrics,  we obtain:

\begin{teo}\label{glavna2}
The magnetic sub-Riemannian geodesic flows of the  metrics $\{ds^2_{sR,a}\}$
 are completely integrable in the non-commutative sense.
The complete algebra of first integrals is generated by the components of the momentum mapping $\Phi_\eta^{ij}=\langle \Phi_\eta,\mathbf E_i\wedge \mathbf E_j\rangle$ and the Hamiltonian function $H_{sR,a}$. For $n\ge 4$, the dimension of invariant tori is equal to $3$.
For $n=3$, the dimension of invariant tori is equal to $2$.
\end{teo}

Note that for $\eta=0$ we get the  integrability of normal sub-Riemannian geodesic flows of the metrics $\{ds^2_{sR,a}\}$.
 Vice verse, the magnetic systems, via reduction, naturally appear in the study of normal sub-Riemannian geodesic flows (e.g., see \cite{BDM, PS}).

A detail study of sub-Riemannian geodesic flows on $SO(3)\cong V_{3,2}$ is given in \cite{BS2016}.
Other related extremal problems on Stiefel varieties are given in \cite{Ju}.
It should also be noted that there is a subtle relationship between elastic problems (as part of sub-Riemannian geometry) and the Kowalevski top (see \cite{Ju2016}).

\begin{rem}[Sub-Riemannian Manakov metrics]\label{primedba*}
The integrability of the flows for Hamiltonian functions $H_{\mathcal H,a}$ (for $a_4=0$) and $H_{\mathcal D_0,a}$, without the magnetic field, also follows from the results of \cite{DGJ2009, DGJ2015, Myk}.
Namely, for the Manakov metrics $\{ds^2_{A,B}\}$ (see Remark \ref{primedba}), the limit $a_1+a_2\to 2a_3$ is equivalent to the
the limit $\beta_2\to\beta_1$, while the limit $a_2\to 0$ is equivalent to the limit $\beta_3\to\beta_2$
(see \eqref{parametri-manakov}).
In the limit $\beta_2\to\beta_1$, the scalar product \eqref{manakov} becomes singular at $\mathfrak v_{1,2}$. However, it is well-defined at
$\mathcal H\vert_o=\mathfrak v_{1,3}\oplus\mathfrak v_{2,3}$,
inducing $SO(n)$--invariant sub-Reimannian metric on the contact distribution $\mathcal H$.
Thus, we get the subclass of sub-Riemannian metrics $\{ds^2_{\mathcal H,a}\,\vert\, a_1+a_2=2a_3, a_4=0\}$.
Similarly,  in the limit $\beta_3\to\beta_2$,
the scalar product \eqref{manakov} becomes singular at $\mathfrak v_{2,3}$. However, it is well-defined at
$\mathfrak d_0=\mathfrak v_{1,2}\oplus\mathfrak v_{1,3}$,
inducing $SO(n)$--invariant sub-Riemannian metric $ds^2_{\mathcal D_0,a}$,
where
\[
a_1=\frac{\beta_1-\beta_2}{\alpha_1-\alpha_3}, \quad a_3=\frac12\frac{\beta_1-\beta_2}{\alpha_1-\alpha_3}-2\frac{\beta_1-\beta_2}{\alpha_1-\alpha_2},
\quad \frac{\beta_1-\beta_2}{\alpha_1-\alpha_3}>0, \quad  \frac{\beta_1-\beta_2}{\alpha_1-\alpha_2}>0.
\]
On the other hand, the Manakov integrals do not depend on the matrix $B$ (see \cite{Ma} for generic matrixes $A$ and \cite{DGJ2009} for matrixes $A$ with multiple eigenvalues), which implies integrability of sub-Riemannian geodesic flows both on $SO(n)$, as well as on the homogeneous spaces of the Lie group $SO(n)$ (see \cite{DGJ2009, BAB}).
\end{rem}

Finally, note that on the compact Lie groups we have a natural construction
of left-invariant (or equivalently, right-invariant) sub-Riemannian structures with completely integrable flows (e.g.,
using chains of subalgebras, see \cite{JSV2024}).  Still, if a sub-Riemannian structure is defined by restricting a left-invariant Riemannian metric to a right-invariant bracket generating distribution,
the construction of integrable examples is a rather difficult task.
One of the first steps in this direction was recently given by Pavlovi\'c and \v Sukilovi\'c in \cite{PS}.

\section{Magnetic pendulum systems on $V_{n,2}$}\label{sec4}

\subsection{Natural mechanical systems on $V_{n,2}$}

We consider natural mechanical systems $(V_{n,2}, L_{\kappa})$ with the kinetic energy given by $SO(n)\times SO(2)$ invariant metrics $ds^2_{1,\kappa}$ and the
potential function $V(\e_1,\e_2)$. According to Lemma \ref{lezandr}, the Lagrangian of the system is
\[
L_\kappa(\e_1,\e_2,\dot \e_1,\dot \e_2)=\frac12\langle \dot \e_1,\dot
\e_1\rangle+\frac12\langle \dot \e_2,\dot
\e_2\rangle+\frac{1+2\kappa}{2+2\kappa}\langle \e_1,\dot
\e_2\rangle\langle \e_2,\dot \e_1\rangle-V(\e_1,\e_2),
\]
while the Hamiltonain is given by
\begin{equation}
H_{\kappa}=\frac12\langle \p_1,\p_1\rangle +\frac12\langle
\p_2,\p_2\rangle -(1+2\kappa)\langle \p_1,\e_2\rangle \langle
\p_2,\e_1\rangle+V(\e_1,\e_2). \label{Ham}
\end{equation}

 Again, by applying equations \eqref{H*} we get:

\begin{lem}\label{prirodni-tok}
The Hamiltonian equations of the Hamiltonian \eqref{Ham} with respect to the twisted symplectic structure
\eqref{twisted} are given by
\begin{equation}\label{MagFlow}
\begin{aligned}
&\dot \e_1=\p_1-(1+2\kappa)\langle \e_1,\p_2\rangle  \e_2, \\
& \dot \e_2=\p_2-(1+2\kappa)\langle \e_2,\p_1\rangle \e_1, \\
& \dot \p_1=-\frac{\partial V}{\partial \e_1} + (1+2\kappa)\langle \p_1,\e_2\rangle \p_2+\eta\,(\p_2-(1+2\kappa)\langle \e_2,\p_1\rangle \e_1)+\lambda_{11} \e_1+ \lambda_{12} \e_2, \\
& \dot \p_2=-\frac{\partial V}{\partial \e_2} + (1+2\kappa)\langle \p_2,\e_1\rangle \p_1-\eta
\, (\p_1-(1+2\kappa)\langle \e_1,\p_2\rangle  \e_2)+\lambda_{12} \e_1+
\lambda_{22} \e_2,
\end{aligned}
\end{equation}
where the Lagrange multipliers $\lambda_{ij}$ are
\begin{align}
\nonumber & \lambda_{11}=-\langle \p_1,\p_1\rangle+\langle
\e_1,\frac{\partial V}{\partial \e_1}\rangle-\eta \langle \p_2,\e_1
\rangle +\eta (1+2\kappa)\langle \e_2,\p_1\rangle,\\
\label{mnozioci} & \lambda_{12}=-\langle \p_1,\p_2\rangle+
\frac12\big(\langle \e_1,\frac{\partial V}{\partial \e_2}\rangle+\langle \e_2,\frac{\partial V}{\partial \e_1}\rangle\big), \\
\nonumber &\lambda_{22}=-\langle \p_2,\p_2\rangle+\langle
\e_2,\frac{\partial V}{\partial \e_2}\rangle+\eta \langle \p_1,\e_2\rangle -\eta (1+2\kappa)\langle \p_2,\e_1\rangle.
\end{align}
\end{lem}
%Note that for the case of the normal metric ($\kappa=0$), the Lagrange multipliers do not depend on the magnetic field.

 In \cite{FeJo_JPA_2012} we have studied
systems with two types of pendulum potentials that we recall below.

\subsection{Pendulum system induced from $G^+_{n,2}$}\label{sec4I}

A class of integrable magnetic pendulum systems on adjoint orbit of compact Lie groups with respect to the magnetic term induced by the
Kostant-Kirilov symplectic form is obtained in \cite{BJ2008}. For the case of the Grassmannian manifold $G^+_{n,2}$ of oriented two-planes in $\R^n$, the potential is given by
\[
v(\e_1\wedge e_2)=\langle \e_1\wedge \e_2,\Xi\rangle,
\]
where $\Xi\in so(n)$. That is why it is a natural problem to consider the pull-back of $v$ with respect to the
submersion \eqref{submersion}:
\begin{equation}\label{pendulumI}
V_I(\e_1,\e_2)=v(\pi(\e_1,e_2))= \langle \e_1\wedge \e_2,\Xi\rangle=-\langle \Xi \e_1,\e_2\rangle.
\end{equation}

The equations are given by \eqref{MagFlow}, where the derivatives of $V_I$ are:
\begin{equation*}\label{LagTop}
\frac{\partial V_I}{\partial \e_1}=\Xi \e_2, \qquad \frac{\partial V_I}{\partial \e_2}=-\Xi \e_1.
\end{equation*}

Since the potential \eqref{pendulumI} is right $SO(2)$--symmetric, the momentum mapping $\Psi$  is the first integral of the system.
The following theorem is proved in \cite{FeJo_JPA_2012} (set $\Xi\mapsto -\Xi$ in Theorem 7) for the case $\eta=0$.
Using the relations
\[
\frac{d}{dt}\big(\e_1 \wedge \e_2\big)=[\Phi_\eta,\e_1\wedge \e_2], \qquad
\dot\Phi_\eta=[\Xi,\e_1 \wedge \e_2],
\]
the proof can be easily modified for the case $\eta\ne 0$.

\begin{teo} \label{Gn2-klatno}
(i) The equations \eqref{MagFlow} with the potential \eqref{pendulumI} imply the following matrix equation with
a spectral parameter $\lambda$
\begin{equation}\label{lax1}
\begin{aligned}
&\frac{d}{dt}
\mathcal L(\lambda)=[\mathcal A(\lambda),\mathcal L(\lambda)],\\
& \mathcal L(\lambda)=\lambda \Phi_\eta + \sqrt{-1}\, ( \e_1\wedge \e_2-\lambda^2
\Xi), \quad \mathcal A(\lambda)=\Phi_\eta-\sqrt{-1}\, \lambda \Xi.
\end{aligned}
\end{equation}

(ii) For a generic $\Xi\in so(n)$, the
magnetic system \eqref{MagFlow} is Liuville integrable. The Noether integral $\Psi$ and the invariants of the $\mathcal L(\lambda)$ form a complete set of $(2n-3)$
involutive polynomials on $(T^*V_{n,2},\omega_\eta)$.
\end{teo}

\subsection{Pendulum system induced from $S^{n-1}\times S^{n-1}$}\label{sec4II}

Let $\gamma_1,\gamma_2\in\R^n$ be two arbitrary fixed unit vectors.
We consider the pendulum-type potential on $V_{n,2}$:
\begin{equation}\label{pendulumII}
V_{II}(\e_1,\e_2)=\chi_1\langle \gamma_1,\e_1\rangle+\chi_2\langle\gamma_2,\e_2\rangle,
\end{equation}
For the case of the Eucledian metric ($\kappa=-1/2$),
the problem represents two independent pendulums systems on $S^{n-1}(\e_1)$ and
$S^{n-1}(\e_2)$ coupled by the holonomic constraint $\langle \e_1,\e_2\rangle=0$.

Now, the potential is not right $SO(2)$--invariant and the equations are given by \eqref{MagFlow}, where the derivatives of $V_{II}$ are
${\partial V_{II}}/{\partial \e_1}=\chi_1 \gamma_1$ and ${\partial V_{II}}/{\partial \e_2}=\chi_2\gamma_2$.

Let us set
\[
\Gamma=(\chi_1\gamma_1,\chi_2\gamma_2), \qquad
\Psi_\eta= \begin{pmatrix}
 0& \Psi-\eta\\
 \eta-\Psi & 0
\end{pmatrix}
\]
Using the relations
\begin{align*}
&\dot X=\Phi_0 X+\kappa X\Psi_0=\Phi_\eta X+\kappa X \Psi_\eta,\\
&\dot \Phi_\eta=X\Gamma^T-\Gamma X^T=\chi_1\e_1\wedge\gamma_1+\chi_2 \e_2\wedge \gamma_2,\\
&\dot \Psi_\eta=\Gamma^T X-X^T \Gamma,
\end{align*}
we get the following Lax representation of the problem, established in \cite{FeJo_JPA_2012} for $\eta=0$.

\begin{teo} \label{Vn2-klatno}  The
equations \eqref{MagFlow}, for the pendulum-type potential \eqref{pendulumII} and $\kappa=1$ imply the
matrix equations
\begin{equation}\label{LASPmag}
\frac{d}{dt} \mathcal L(\lambda)=[\mathcal L(\lambda),\mathcal A(\lambda)]
\end{equation}
with a spectral parameter $\lambda$, where $\mathcal L(\lambda)$ and
$\mathcal A(\lambda)$ are $so(n+2)$ matrices given by:
\begin{align*}
& \mathcal L(\lambda)=\begin{pmatrix}
-\lambda\Phi_\eta & \e_1+\lambda^2\chi_1\gamma_1 & \e_2+\lambda^2\chi_2\gamma_2 \\
-\e_1^T-\lambda^2\chi_1\gamma_1^T & 0& \lambda (\Psi-\eta)\\
-\e_2^T-\lambda^2\chi_2\gamma_2^T & \lambda(\eta-\Psi) & 0
\end{pmatrix}, \\
& \mathcal A(\lambda)=\begin{pmatrix}
-\Phi_\eta & \lambda\chi_1\gamma_1 & \lambda\chi_2\gamma_2 \\
-\lambda\chi_1\gamma_1^T & 0& \Psi-\eta\\
-\lambda\chi_2\gamma_2^T &\eta-\Psi & 0
\end{pmatrix}.
\end{align*}
\end{teo}

From the Lax representation \eqref{LASPmag}, we obtain:

\begin{teo}\label{main}
The magnetic pendulum system \eqref{MagFlow}, \eqref{pendulumII} for $\kappa=1$ is completely
integrable in the noncommutative sense.
The dimension of invariant tori $\delta$ is given by:
\begin{eqnarray*}
 &                   & \delta=5,      \qquad\,\,\,\,            n\ge 5,\\
&\gamma_1=\gamma_2\qquad\,\,\,\,& \delta=4, \qquad\,\,\,\,  n=4,\\
  &                  & \delta=3, \qquad\,\,\,\,  n=3.
\end{eqnarray*}
\begin{eqnarray*}
   &                         & \delta=7, \qquad\,\,\,\, n\ge 6,\\
&\gamma_1\ne\gamma_2\qquad\,\,\,\,& \delta=6,\qquad\,\,\,\,   n=5,\\
 &                           & \delta=5,\qquad\,\,\,\,   n=4,\\
  &                          & \delta=3,\qquad\,\,\,\,   n=3.
\end{eqnarray*}
\end{teo}

We will prove Theorem \ref{main} in a separate paper.

\subsection{$n=3$ and the heavy top with a gyrostat}\label{sec4III}

In the case $n=3$, we write the matrix $\Xi\in so(3)$ in the form
\begin{equation}\label{Xi}
\Xi=
\chi_3 \begin{pmatrix}
0 & -\gamma_3 & \gamma_2 \\
\gamma_3 & 0 & -\gamma_1 \\
-\gamma_2 & \gamma_1 & 0
\end{pmatrix}, \qquad \gamma_1^2+\gamma_2^2+\gamma_3^2=1,
\end{equation}
and set $\gamma=(\gamma_1,\gamma_2,\gamma_3)$. Then
\[
V_I=-\langle \Xi \e_1,\e_2\rangle=-\chi_3 \langle \gamma\times \e_1,\e_2\rangle=-\chi_3 \langle \gamma,\e_2\times \e_1\rangle=\chi_3 \langle\gamma,\mathbf \e_3\rangle.
\]

Let $\Gamma$ be the vector $\gamma$ in the moving reference frame and $\chi=(0,0,\chi_3)$. From \eqref{rigid-body}, the Hamiltonian of the system in the moving frame reads
\[
H_{I,\kappa}=\frac12(M_1^2+M_2^2+(1+\kappa)M_3^2)+\langle \Gamma,\chi\rangle,
\]
and we get:

\begin{prp}\label{stav2}
For $n=3$, the magnetic flows \eqref{MagFlow} with respect to potential \eqref{pendulumI}, \eqref{Xi} represents a motion of a system Lagrange top+gyrostat  around a fixed point in the presence of the homogeneous gravitational field.  The inertia operator
of the body is $I=\diag(1,1,1/(1+\kappa))$, the gyrostat momentum is $\mathbf L=(0,0,-\eta)$, and the position of the center of mass of the system body+gyrostat, multiplied by the mass of the system $m$  and the gravitational constant $g$,
 is given by  $\chi=(0,0,\chi_3)$.
\end{prp}

For $n=3$, the Lax representation \eqref{lax1} of the Lagrange top with a gyrostat is the Lax representation with respect to the fixed coordinate system.
In particular, for $\eta=0$, this is a dual Lax representation to the Lax representation of the Lagrange top obtained in \cite{RM}.

Now we consider the system with the pendulum-type potential \eqref{pendulumII}.
Assume $\gamma_1=\gamma_2=\gamma$, and let $\Gamma$ be the vector $\gamma$ considered in the moving reference frame. Then we
have
\[
V_{II}=\langle \chi_1\e_1+ \chi_2 \e_2,\gamma \rangle=\langle \chi,\Gamma\rangle, \qquad \chi=(\chi_1,\chi_2,0),
\]
and the Hamiltonian of the system in the moving frame reads
\[
H_{II,\kappa}=\frac12(M_1^2+M_2^2+(1+\kappa)M_3^2)+\langle \Gamma,\chi\rangle,
\]
The inertia operator \eqref{IO} for $\kappa=1$ is the inertia operator of the Kowalevski top. Thus, we get:

\begin{prp}\label{stav3}
For $n=3$ and $\kappa=1$, the magnetic flows \eqref{MagFlow} with respect to potential \eqref{pendulumII} with $\gamma_1=\gamma_2$ represents a motion of a system Kowalevski top top+gyrostat  around a fixed point in the presence of the homogeneous gravitational field.  The inertia operator
of the body is $I=\diag(1,1,1/2)$, the gyrostat momentum is $\mathbf L=(0,0,-\eta)$,
and the position of the center of mass of the system body+gyrostat, multiplied by the mass of the system $m$  and the gravitational constant $g$,
 is given by $\chi=(\chi_1,\chi_2,0)$.
\end{prp}

Thus, for $n=3$, the Lax representation  \eqref{LASPmag} for
the Kowalewski gyrostat is a dual Lax representation with respect to those given in  \cite{RST, D1997, BRS, ST}.

If $\gamma_1\ne\gamma_2$, the Lax representation \eqref{LASPmag} implies the integrability of so called Kowalewski top with a gyrostat under the influence of two
homogeneous force fields.

\subsection*{Acknowledgments}
The author is very grateful to the referees for useful remarks.
This research was supported
by the Serbian Ministry of Science, Technological Development and
Innovation through Mathematical Institute of Serbian Academy of Sciences
and Arts and it is a part of the proposal IntegraRS of the Science Fund of
Serbia.

\end{document}